\documentclass[twocolumn,10pt,a4]{asme2ej_DETC2018}

\usepackage{amsfonts}
\usepackage{graphicx}
\usepackage{epstopdf}
\usepackage{amsmath, amssymb, bm}
\usepackage[fixamsmath,disallowspaces]{mathtools}

\newtheorem{remark}{Remark}
\newtheorem{definition}{Definition}

\begin{document}

\parindent 0pt \setcounter{topnumber}{9} \setcounter{bottomnumber}{9} %
\renewcommand{\textfraction}{0.00001}

\renewcommand {\floatpagefraction}{0.999} \renewcommand{\textfraction}{0.01} %
\renewcommand{\topfraction}{0.999} \renewcommand{\bottomfraction}{0.99} %
\renewcommand{\floatpagefraction}{0.99} \setcounter{totalnumber}{9}

\title{ 
\vspace{-5ex}
A Screw Approach to the Approximation of the Local Geometry of the Configuration Space and of the set of Configurations of Certain Rank of Lower Pair Linkages
\vspace{-8mm}
} 
\author{
	Andreas M\"uller
	\affiliation{Johannes Kepler University, Linz, Austria, a.mueller@jku.at}
\vspace{-5ex}
} 

\maketitle 

\begin{abstract}
\textit{Abstract--} A motion of a mechanism is a curve in its configuration
space (c-space). Singularities of the c-space are kinematic singularities of
the mechanism. Any mobility analysis of a particular mechanism amounts to
investigating the c-space geometry at a given configuration. A higher-order
analysis is necessary to determine the finite mobility. To this end, past
research lead to approaches using higher-order time derivatives of loop
closure constraints assuming (implicitly) that all possible motions are
smooth. This continuity assumption limits the generality of these methods.
In this paper an approach to the higher-order local mobility analysis of
lower pair multi-loop linkages is presented. This is based on a higher-order
Taylor series expansion of the geometric constraint mapping, for which a
recursive algebraic expression in terms of joint screws is presented. An
exhaustive local analysis includes analysis of the set of constraint
singularities (configurations where the constraint Jacobian has certain
corank). A local approximation of the set of configurations with certain
rank is presented, along with an explicit expression for the differentials
of Jacobian minors in terms of instantaneous joint screws. The c-space and
the set of points of certain corank are therewith locally approximated by an
algebraic variety determined algebraically from the mechanism's screw system.

Results are shown for a simple planar 4-bar linkage, which exhibits a
bifurcation singularity, and for a planar three-loop linkage exhibiting a
cusp in c-space. The latter cannot be treated by the higher-order local
analysis methods proposed in the literature.


\vspace{-3ex}

\end{abstract} 

\newpage%

\section{Introduction}

Mobility and singularity analysis of linkages remains a topic of extensive
research. Ideally, such an analysis reveals a global picture of the possible
motions and the existence of singularities as well as the motion
characteristics when passing through singularities. Yet, a universal
approach to this problem is not available, and various approaches where
developed that applicable for specific types of problems. On the other hand,
screw and Lie group theory allow for a kinematics modeling that is
advantageous for a general higher-order \emph{local} mobility and
singularity analysis. Research on the \emph{global} analysis (of mechanisms
comprising algebraic joints) also progressed recently using methods from
computational algebraic geometry \cite{HustyRobotica2007}. In particular the
Linear Implicitization Algorithm proposed in \cite{WalterHusty2010} is a
method to derive algebraic constraint equations in a systematic way and in a
form beneficial for the computational analysis. Using local joint frames
according to the Denavit-Hartenberg parameterization, the algebraic loop
constraints are derived by substitution of the trigonometric relations as
shown in \cite{KongJMR2017}. Another approach in terms of DH parameters was
presented in \cite{LiSchichoMMT2015}. While computational algebraic methods
allow for global mobility analysis of certain mechanisms, the application of
algebraic methods is limited due to the complexity of the involved
algorithms. A local analysis, on the other hand, is deemed generally
applicable. Various approaches to the higher-order local analysis, aim to
approximate possible finite motion through a given (general) configuration,
were reported in \cite%
{Bustos2012,Gallardo2001,Lerbet1999,Rico1999,JMR2016,JMR2018}. Although
involving different formulations of loop constraints, all these methods rest
on using higher-order derivatives of the geometric loop constraints. The
basic idea is to identify tangents to finite curves in configuration space
(c-space), which then indicate possible motions and the local finite degree
of freedom (DOF), that locally approximate the c-space. Simply speaking,
these methods determine generalized \emph{velocity} vectors for which there
exist a generalized acceleration, jerk, jounce, etc. satisfying the
constraints of order two, three, four, etc. This implies that the
so-determined vectors are tangents to \emph{smooth} curves in c-space.
Non-smooth curves are hence excluded. Therefore these tangents (forming a
cone rather than a vector space) are said to form the \emph{kinematic
tangent cone} in order to emphasize that they correspond to finite motions
with non-zero velocity. Given the fact that the c-space singularities of
virtually all mechanisms that were reported in the literature do allow for
smooth motions between any two points, e.g. bifurcations, it may not be
surprising that the above approaches to a higher-order local analysis are
often assumed to provide conclusive statements. Until recently, the only
mechanisms exhibiting a c-space singularity that does not allow for smooth
motions, namely a cusp, through this singularity was the double-Watt linkage
reported in \cite{ConnellyServatius1994}. Other spatial linkages were
reported in \cite{PabloMMT2019}. Such linkages are said to be $k$th-order
rigid while being finitely mobile \cite{ConnellyServatius1994}. Clearly,
other methods for the higher-order local analysis are needed to account for
such situations. Such a method is the local approximation of the c-space
independently of possible finite motions. This should seem to be the most
natural approach, as finite motions of a linkage are encoded in the c-space
geometry, but the resulting equations can be difficult to use. This is
addressed in the present paper. A higher-order series expansion of the
constraints is used to this end. Relevant prior publications are \cite%
{Chen2011,Karsai2001,Milenkovic2012}. In \cite{Chen2011} a second-order
analysis was presented using the series expansion of the constraints. In 
\cite{Karsai2001} and \cite{Milenkovic2012} higher-order series expansions
were reported and applied to solve inverse kinematic problems.

The original contribution of this paper is a method for the local
approximation of the geometry of the c-space in section \ref{secLocCSpace}
and the recursive relations for the differentials in the series expansion of
the constraint mapping in terms of joint screws, in section \ref{secDiff}.
An approach to the approximation of the subvariety of points of certain rank
is introduced in section \ref{secLocLk}. A closed form of the required
differentials of the Jacobian minors is presented in section \ref{secDiff}.
Together with the prior publications \cite{JMR2016,JMR2018}, this paper
completes the conceptual and algorithmic framework for an exhaustive local
analysis of general multi-loop linkages in terms of joint screws. The
constraint formulation and definition of kinematic singularities, constraint
singularities, and c-space singularities are recalled in section \ref%
{secConstraints}, and that of the kinematic tangent cone in section \ref%
{secKinTangCone}. Two examples are shown in section \ref{secExamples}. The
first example is the simple planar 4-bar linkage, which only serves to show
that the method reproduces the known result. The second example is the
planar three-loop double-Watt linkage reported in \cite%
{ConnellyServatius1994,JMR2018}, which exhibits a cusp singularity in c-space%
\vspace{-2ex}%

\section{Constraints, Configuration Space, Singularities%
\label{secConstraints}%
}

First consider a kinematic loop comprising $n$ 1-DOF lower pair joints. The
joint variable $q_{i}$ is assigned to joint $i$, where $q_{i}\in \mathbb{T}$
if the joint is a revolute or screw joint, and $q_{i}\in \mathbb{R}$ for a
prismatic joint. The vector of joint variables is denoted with $\mathbf{q}%
\in {\mathbb{V}}^{n}={\mathbb{T}}^{n_{\text{R}}}\times {\mathbb{R}}^{n_{%
\text{P}}}$, where $n_{\text{R}}$ is the number of revolute joints and $n_{%
\text{P}}$ the number of prismatic and helical joints. Denote with $\mathbf{Y%
}_{i}$ the screw coordinates associated to joint $i$, expressed in the world
frame, determined in the zero reference configuration $\mathbf{q}=\mathbf{0}$%
. Presuming that $\mathbf{0}\in V$, the geometric loop closure constraints
can then be expressed as $f\left( \mathbf{q}\right) =\mathbf{I}$ where the
constraint mapping is%
\begin{equation}
f\left( \mathbf{q}\right) =\exp \left( \mathbf{Y}_{1}q_{1}\right) \exp
\left( \mathbf{Y}_{2}q_{2}\right) \cdot \ldots \cdot \exp \left( \mathbf{Y}%
_{n}q_{n}\right) .  \label{f}
\end{equation}%
This formula is referred to as the product of exponentials (POE). It is
usually attributed to Brockett \cite{Brockett1984}, but was already
presented by Herve in \cite{Herve1982}.

\begin{remark}
\label{rem1}%
The joint screw coordinates $\mathbf{Y}_{i}$ in (\ref{f}) are determined in
an assembly configuration of the linkage, so that $\mathbf{q}=\mathbf{0}$ is
an admissible configuration. Moreover, w.l.o.g. throughout the paper, the
zero reference configuration is chosen as the configuration of interest.
\end{remark}

The velocity loop constraints are%
\begin{equation}
\mathbf{0}=\mathbf{J}\left( \mathbf{q}\right) \dot{\mathbf{q}}=\sum_{i\leq n}%
\mathbf{S}_{i}\left( \mathbf{q}\right) \dot{q}_{i}  \label{VelConstr}
\end{equation}%
where the instantaneous screw coordinate vector of joint $i$ 
\begin{equation}
\mathbf{S}_{i}%
\hspace{-0.5ex}%
\left( \mathbf{q}\right) :=\mathbf{Ad}_{g_{i}\left( \mathbf{q}\right) }%
\mathbf{Y}_{i}  \label{Si}
\end{equation}%
with $g_{i}%
\hspace{-0.5ex}%
\left( \mathbf{q}\right) :=\exp (\mathbf{Y}_{1}q_{1})\cdot \ldots \cdot \exp
(\mathbf{Y}_{i}q_{i})$. The $\mathbf{S}_{i}$ are the columns of the
geometric constraint Jacobian.

Now consider a general linkage with multiple kinematic loops. This is
modeled as a system of $\gamma $ topologically independent loops --the
fundamental cycles (FCs) of the topological graph \cite{Robotica2017}. For
each of these $\gamma $ loops a system of geometric constraints 
\begin{equation}
f_{l}\left( \mathbf{q}\right) =\mathbf{I},l=1,\ldots ,\gamma  \label{fl}
\end{equation}%
is introduced. The corresponding velocity constraints 
\begin{equation}
\mathbf{J}_{l}\left( \mathbf{q}\right) \dot{\mathbf{q}}=\mathbf{0}%
,l=1,\ldots ,\gamma  \label{VelConstr2}
\end{equation}%
are summarized in the form (\ref{VelConstr}) using the overall Jacobian%
\begin{equation}
\mathbf{J}=\left( 
\begin{array}{c}
\mathbf{J}_{1} \\ 
\vdots \\ 
\mathbf{J}_{\gamma }%
\end{array}%
\right) .  \label{J}
\end{equation}

The c-space is defined by geometric loop constraints as%
\begin{equation}
V:=\left\{ \mathbf{q}\in {\mathbb{V}}^{n}|f_{l}\left( \mathbf{q}\right) =%
\mathbf{I},l=1,\ldots ,\gamma \right\} .  \label{V}
\end{equation}%
The \emph{local DOF} of the linkage at $\mathbf{q}\in V$, denoted with $%
\delta _{\mathrm{loc}}\left( \mathbf{q}\right) $, is the local dimension of
the solution set $V$ at $\mathbf{q}$. The \emph{differential DOF} of the
linkage is $\delta _{\text{diff}}\left( \mathbf{q}\right) :={n-\mathrm{rank}~%
}\mathbf{J}\left( \mathbf{q}\right) $. The configuration $\mathbf{q}$ is
said to have rank $r$ when $r={\mathrm{rank}~}\mathbf{J}\left( \mathbf{q}%
\right) $.

A configuration $\mathbf{q}\in {\mathbb{V}}^{n}$ is a \emph{constraint
singularity} iff ${\mathrm{rank}~}\mathbf{J}\left( \mathbf{q}\right)
<d_{\max }$, where $d_{\max }$ is the maximal rank of $\mathbf{J}$ for $%
\mathbf{q}\in {\mathbb{V}}^{n}$ (of the open kinematic chain). It is a \emph{%
kinematic singularity} iff there is no neighborhood of $\mathbf{q}$ in $V$
where $\delta _{\text{diff}}$ is constant. It is \emph{c-space singularity}
iff $V$ is not a smooth manifold at $\mathbf{q}$. A c-space singularity is
always a constraint singularity and a kinematic singularity. A kinematic
singularity is always a constraint singularity but may not be a c-space
singularity. A constraint singularity may not be a c-space nor a kinematic
singularity.

\begin{remark}
An upper bound for the maximal rank of $\mathbf{J}$ is the dimension of the
screw algebra generated by the joint screws (i.e. the involutive closure of
the screw system). This is the vector space $\mathrm{span}~\left( \mathbf{Y}%
_{i},[\mathbf{Y}_{i},\mathbf{Y}_{j}],[\mathbf{Y}_{i},[\mathbf{Y}_{j},\mathbf{%
Y}_{l}]\right) $, which is easily determined by nested Lie brackets/screw
products of joint screws $\mathbf{Y}_{i}$.
\end{remark}

\section{%
\label{secKinTangCone}%
Kinematic Tangent Cone revisited}

\subsection{Kinematic Tangent Cone to Smooth Curves}

The basic idea of a kinematic tangent cone was presented in \cite{Lerbet1999}
and then formalized in \cite{JMR2016} as a means to identify the tangent
properties of smooth curves through a given point in c-space. As discussed
in \cite{JMR2018}, this reveals the local c-space geometry at regular and
singular points of $V$, except at cusps, however.

\begin{definition}
Let $V$ be the c-space of a mechanism. The \emph{kinematic tangent cone} of $%
V$ at $\mathbf{q}\in V$, denoted $C_{\mathbf{q}}^{\mathrm{K}}$, is the set
of tangent vectors to smooth arcs in $V$ passing through $\mathbf{q}$.
\end{definition}

The point of departure for constructing the kinematic tangent cone is the
system of velocity loop constraints (\ref{VelConstr2}). To this end, the
definition%
\begin{equation}
H_{l}^{\left( 1\right) }%
\hspace{-0.6ex}%
\left( \mathbf{q},\dot{\mathbf{q}}\right) :=\mathbf{J}_{l}\left( \mathbf{q}%
\right) \dot{\mathbf{q}}
\end{equation}%
will be helpful. Solutions $\mathbf{x}\in {\mathbb{R}}^{n}$ of $\mathbf{0}%
=H_{l}^{\left( 1\right) }%
\hspace{-0.6ex}%
\left( \mathbf{q},\mathbf{x}\right) ,l=1,\ldots ,\gamma $ are admissible
1st-order infinitesimal joint motions. A $\dot{\mathbf{q}}$ is a tangent
vector to a smooth curve through $\mathbf{q}\in V$ iff there are
corresponding $\ddot{\mathbf{q}},\dddot{\mathbf{q}},\ldots $ that satisfy
all higher-order constraints%
\begin{equation}
H_{l}^{\left( i\right) }%
\hspace{-0.5ex}%
(\mathbf{q},\dot{\mathbf{q}},\ldots ,\mathbf{q}^{\left( i\right) })=\mathbf{0%
}  \label{Hi}
\end{equation}%
with $H_{l}^{\left( i\right) }%
\hspace{-0.5ex}%
(\mathbf{q},\dot{\mathbf{q}},\ldots ,\mathbf{q}^{\left( i\right) }){:=}\frac{%
d^{i-1}}{dt^{i-1}}H_{l}^{\left( 1\right) }%
\hspace{-0.6ex}%
\left( \mathbf{q},\dot{\mathbf{q}}\right) $, using the notation $\mathbf{q}%
^{(i)}%
\hspace{-0.5ex}%
:=%
\hspace{-0.2ex}%
\frac{d^{i}}{dt^{i}}\mathbf{q}$. The kinematic tangent cone, $C_{\mathbf{q}%
}^{\mathrm{K}}$, is then determined by the sequence \cite{JMR2016}%
\vspace{-1ex}%
\begin{equation}
{C_{\mathbf{q}}^{\text{K}}V}=K_{\mathbf{q}}^{\kappa }\subset \ldots \subset
K_{\mathbf{q}}^{3}\subset K_{\mathbf{q}}^{2}\subset {K_{\mathbf{q}}^{1}}%
\vspace{-1ex}
\label{CqV}
\end{equation}%
where the $i$th-order cone $K_{\mathbf{q}}^{i}$ is defined as%
\vspace{-1ex}%
\begin{equation}
\begin{array}{lll}
K_{\mathbf{q}}^{i}:=\{\mathbf{x}|\exists \mathbf{y},\mathbf{z},\ldots \in {%
\mathbb{R}}^{n}: & H_{l}^{\left( 1\right) }%
\hspace{-0.6ex}%
\left( \mathbf{q},\mathbf{x}\right) =\mathbf{0}, &  \\ 
& H_{l}^{\left( 2\right) }%
\hspace{-0.6ex}%
\left( \mathbf{q},\mathbf{x},\mathbf{y}\right) =\mathbf{0}, &  \\ 
& H_{l}^{\left( 3\right) }%
\hspace{-0.6ex}%
\left( \mathbf{q},\mathbf{x},\mathbf{y},\mathbf{z}\right) =\mathbf{0}, &  \\ 
\multicolumn{1}{c}{} & \multicolumn{1}{c}{\cdots} & \multicolumn{1}{c}{} \\ 
\multicolumn{1}{c}{} & \multicolumn{1}{c}{H_{l}^{\left( i\right) }%
\hspace{-0.6ex}%
\left( \mathbf{q},\mathbf{x},\mathbf{y},\mathbf{z,\ldots }\right) =\mathbf{0}%
,l=1,\ldots ,\gamma} & \multicolumn{1}{c}{\}.}%
\end{array}
\label{Ki}
\end{equation}%
$K_{\mathbf{q}}^{i}$ is the set of admissible $i$th-order motions. The
sequence (\ref{CqV}) terminates with a finite $\kappa $, and $\delta _{\text{%
loc}}\left( \mathbf{q}\right) =K_{\mathbf{q}}^{\kappa }$ is the local (i.e.
finite) DOF at $\mathbf{q}$. Vectors $\mathbf{x}\in K_{\mathbf{q}}^{i}$ for $%
i<\kappa $ are not necessarily tangents to finite motions. Crucial for the
analysis is the fact that the higher-order constraint mappings (\ref{Hi})
can be determined explicitly \cite{MMTDerivatives,Gallardo2001,Rico1999},
and more efficiently recursively \cite{MMTConstraints,JMR2018} without the
need for symbolic manipulation (e.g. with computer algebra systems).

\subsection{Kinematic Tangent Cone to Smooth Curves with constant Rank}

An exhaustive mobility analysis further involves calculating the
differential DOF $\delta _{\text{diff}}\left( \mathbf{q}\right) $ at $%
\mathbf{q}\in V$ and whether this changes in any neighborhood, which
indicates a kinematic singularity. The differential DOF is determined by the
rank of the $6\gamma \times n$ Jacobian $\mathbf{J}\left( \mathbf{q}\right) $
in (\ref{J}). If the rank is constant in a neighborhood of $\mathbf{q}\in V$%
, then $\mathbf{q}$ is a regular configuration of the linkage. But it may be
a constraint singularity so that locally, in nearby regular configurations, $%
\mathbf{J}$ has no full rank. Identification of singularities and
determination of the differential DOF in general configurations thus amounts
to analyze the local geometry of the set of points with ${\mathrm{rank}~}%
\mathbf{J}\left( \mathbf{q}\right) $. While $\delta _{\text{diff}}\left( 
\mathbf{q}\right) $ is readily known from the joint screw system, deciding
whether it is locally constant is difficult. The concept of a tangent cone
can be adopted to analyze curves passing through $\mathbf{q}\in V$ for which 
$\delta _{\text{diff}}$ has at least a certain value. To this end, the $%
\bm{\alpha}%
\bm{\beta}%
$-minor of $\mathbf{J}$ of order $k$ \cite{JMR2018}%
\begin{equation}
m_{%
\bm{\alpha}%
\bm{\beta}%
}(\mathbf{q}):=\det \mathbf{J}_{%
\bm{\alpha}%
\bm{\beta}%
}(\mathbf{q})=\left\vert \mathbf{S}_{%
\bm{\alpha}%
\beta _{1}}(\mathbf{q})\cdots \mathbf{S}_{%
\bm{\alpha}%
\beta _{k}}(\mathbf{q})\right\vert
\end{equation}%
is introduced as the determinant of the matrix $\mathbf{J}_{%
\bm{\alpha}%
\bm{\beta}%
}$ constructed using the rows and columns of $\mathbf{J}$ according to the
row and column indexes $%
\bm{\alpha}%
=\{\alpha _{1},\ldots ,\alpha _{k}\}$ and $%
\bm{\beta}%
=\{\beta _{1},\ldots ,\beta _{k}\}$, with $\alpha _{i}\in \{1,\ldots
,6\gamma \},\beta _{j}\in \{1,\ldots ,n\}$. A point $\mathbf{q}\in {\mathbb{V%
}}^{n}$ is of rank less than $k$ iff $m_{%
\bm{\alpha}%
\bm{\beta}%
}(\mathbf{q})=0$ for any $k$-minor, i.e. $\left\vert 
\bm{\alpha}%
\right\vert 
\hspace{-0.5ex}%
=%
\hspace{-0.5ex}%
\left\vert 
\bm{\beta}%
\right\vert 
\hspace{-0.5ex}%
=%
\hspace{-0.3ex}%
k$. The variety of points of rank less than $k$ is thus 
\begin{align}
L_{k}\,{=}& \left\{ \mathbf{q}\in {\mathbb{V}}^{n}|f_{l}\left( \mathbf{q}%
\right) =\mathbf{I},m_{%
\bm{\alpha}%
\bm{\beta}%
}(\mathbf{q})=0,\right.  \label{Lk} \\
& \ \ \ \ \ \ \ \ \ \ \ \ \ \ 
\bm{\alpha}%
\subseteq \{1,\ldots ,6\gamma \},%
\bm{\beta}%
\subseteq \{1,\ldots ,n\},\left\vert 
\bm{\alpha}%
\right\vert 
\hspace{-0.5ex}%
=%
\hspace{-0.5ex}%
\left\vert 
\bm{\beta}%
\right\vert 
\hspace{-0.5ex}%
=%
\hspace{-0.3ex}%
k,  \notag \\
& \ \ \ \ \ \ \ \ \ \ \ \ \ \ \left. l=1,\ldots ,\gamma \right\} \subseteq V%
\vspace{-2ex}%
.  \notag
\end{align}%
For a mechanism with $\gamma $ FC, the number of minors of order $k$ in (\ref%
{Lk}) is $\binom{6\gamma }{k}\cdot \binom{n}{k}$.

\begin{definition}
The kinematic tangent cone to $L_{k}\subseteq V$ at $\mathbf{q}\in V$,
denoted with ${C_{\mathbf{q}}^{\text{K}}}L_{k}$, is the set of tangent
vectors to smooth arcs in $V$ passing through $\mathbf{q}$ at which all
points have rank less than $k$.
\end{definition}

The kinematic tangent cone to $L_{k}$ is determined as%
\begin{equation}
{C_{\mathbf{q}}^{\text{K}}}L_{k}=K_{\mathbf{q}}^{k,\kappa }\subset \ldots
\subset K_{\mathbf{q}}^{k,3}\subset K_{\mathbf{q}}^{k,2}\subset {K_{\mathbf{q%
}}^{k,1}}%
\vspace{-2ex}
\label{CqR}
\end{equation}%
with%
\vspace{-1ex}%
\begin{align}
& K_{\mathbf{q}}^{k,i}:=\{\mathbf{x}|\exists \mathbf{y},\mathbf{z},\ldots
\in {\mathbb{R}}^{n}:  \label{Kki} \\
& \ \ \ \ \ \ \ \ \ \ \ \ \ \ \ \ \ \ 
\begin{array}{l}
\vspace{-3.6ex}
\\ 
H_{l}^{\left( 1\right) }%
\hspace{-0.6ex}%
\left( \mathbf{q},\mathbf{x}\right) =0,M_{%
\bm{\alpha}%
\bm{\beta}%
}^{\left( 1\right) }%
\hspace{-0.6ex}%
\left( \mathbf{q},\mathbf{x}\right) =\mathbf{0},%
\vspace{0.3ex}
\\ 
H_{l}^{\left( 2\right) }%
\hspace{-0.6ex}%
\left( \mathbf{q},\mathbf{x},\mathbf{y}\right) =0,M_{%
\bm{\alpha}%
\bm{\beta}%
}^{\left( 2\right) }%
\hspace{-0.6ex}%
\left( \mathbf{q},\mathbf{x},\mathbf{y}\right) =\mathbf{0},%
\vspace{0.3ex}
\\ 
H_{l}^{\left( 3\right) }%
\hspace{-0.6ex}%
\left( \mathbf{q},\mathbf{x},\mathbf{y},\mathbf{z}\right) =0,M_{%
\bm{\alpha}%
\bm{\beta}%
}^{\left( 3\right) }%
\hspace{-0.6ex}%
\left( \mathbf{q},\mathbf{x},\mathbf{y},\mathbf{z}\right) =\mathbf{0},%
\vspace{0.3ex}
\\ 
\ \ \ \ \ \ \ \ \ \ \ \ \ \ \cdots \\ 
H_{l}^{\left( i\right) }%
\hspace{-0.6ex}%
\left( \mathbf{q},\mathbf{x},\mathbf{y},\mathbf{z,\ldots }\right) =0,M_{%
\bm{\alpha}%
\bm{\beta}%
}^{\left( i\right) }%
\hspace{-0.6ex}%
\left( \mathbf{q},\mathbf{x},\mathbf{y},\mathbf{z,\ldots }\right) =\mathbf{0}%
,%
\vspace{0.3ex}
\\ 
\bm{\alpha}%
\subseteq \{1,\ldots ,6\gamma \},%
\bm{\beta}%
\subseteq \{1,\ldots ,n\},\left\vert 
\bm{\alpha}%
\right\vert 
\hspace{-0.5ex}%
=%
\hspace{-0.5ex}%
\left\vert 
\bm{\beta}%
\right\vert 
\hspace{-0.5ex}%
=k, \\ 
l=1,\ldots ,\gamma \}%
\end{array}%
\vspace{-2ex}%
\vspace{-2ex}
\notag
\end{align}%
where%
\vspace{-2ex}%
\begin{eqnarray}
M_{\mathbf{ab}}^{\left( 1\right) }%
\hspace{-0.6ex}%
\left( \mathbf{q},\dot{\mathbf{q}}\right) {:=} &&\frac{d}{dt}m_{\mathbf{ab}}(%
\mathbf{q})  \label{Mab} \\
M_{\mathbf{ab}}^{\left( 2\right) }%
\hspace{-0.6ex}%
\left( \mathbf{q},\dot{\mathbf{q}},\ddot{\mathbf{q}}\right) {:=} &&\frac{%
d^{2}}{dt^{2}}m_{\mathbf{ab}}(\mathbf{q})  \notag \\
M_{\mathbf{ab}}^{\left( i\right) }%
\hspace{-0.4ex}%
(\mathbf{q},\dot{\mathbf{q}},\ldots ,\mathbf{q}^{\left( i\right) }){:=} &&%
\frac{d^{i}}{dt^{i}}m_{\mathbf{ab}}(\mathbf{q}).%
\vspace{-2ex}
\notag
\end{eqnarray}%
$K_{\mathbf{q}}^{k,i}$ is the set of differential motions of order $i$ where 
${\mathrm{rank}~}\mathbf{J}<k$. These motions are not necessarily tangent to
finite curves. The ${C_{\mathbf{q}}^{\text{K}}}L_{k},1\leq k<{\mathrm{rank}~}%
\mathbf{J}\left( \mathbf{q}\right) $ give rise to a local stratification of $%
V$ according to the rank at smooth curves through $\mathbf{q}\in V$.

For the actual evaluation, it is again important that all derivatives (\ref%
{Mab}) are available recursively or, if desired, in closed form. Denoting
with $\mathbf{a}=\left( a_{1},a_{2},\ldots ,a_{n}\right) \in \mathbb{N}^{n}$
a multi-index, with norm $\left\vert \mathbf{a}\right\vert
:=a_{1}+a_{2}+\ldots +a_{n}$, and $\mathbf{S}_{%
\bm{\alpha}%
\beta _{j}}^{(a_{j})}=\frac{d^{a_{j}}}{dt^{a_{j}}}\mathbf{S}_{%
\bm{\alpha}%
\beta _{j}}$, the derivatives of the minors are given explicitly as\cite%
{JMR2018}%
\begin{equation}
\frac{d^{\nu }m_{%
\bm{\alpha}%
\bm{\beta}%
}}{dt^{\nu }}=\sum_{\left\vert \mathbf{a}_{k}\right\vert =\nu }\left\vert 
\mathbf{S}_{%
\bm{\alpha}%
\beta _{1}}^{\left( a_{1}\right) }\ \mathbf{S}_{%
\bm{\alpha}%
\beta _{2}}^{(a_{2})}\ \cdots \ \mathbf{S}_{%
\bm{\alpha}%
\beta _{k}}^{\left( a_{k}\right) }\right\vert \frac{\nu !}{\mathbf{a}%
_{k}!n_{1}!\cdots n_{\nu }!}  \label{dmdt}
\end{equation}%
where the multi-factorial is defined as $\mathbf{a}!=a_{1}!a_{2}!\cdots
a_{n}!$, and $n_{i}:=|\{a_{j}|a_{j}=i\}|$ is the number of times a
derivative of $\mathbf{S}$ of degree $i$ occurs. The derivatives $\mathbf{S}%
_{%
\bm{\alpha}%
\beta _{i}}^{(a_{i})}$ are given recursively \cite{MMTConstraints,JMR2016}.

\section{%
\label{secLocCSpace}%
Local Approximation of the C-Space}

The c-space is an analytic variety globally defined by the constraint
mapping (\ref{fl}). Replacing the constraint mappings of a multi-loop
linkage by the truncated $\nu $th-order series expansion 
\begin{equation}
f_{l}\left( \mathbf{q}+\mathbf{x}\right) =f_{l}\left( \mathbf{q}\right)
+\sum_{1\leq k\leq \nu }\frac{1}{k!}\mathrm{d}^{k}f_{l,\mathbf{q}}\left( 
\mathbf{x}\right) 
\vspace{-1ex}
\label{fTaylor}
\end{equation}%
gives rise to a local $k$th-order approximation of $V$ at $\mathbf{q}\in V$,
where $\mathrm{d}^{k}f_{l\mathbf{q}}\left( \mathbf{x}\right) $ is the $k$th
differential of $f_{l}$ at $\mathbf{q}\in {\mathbb{V}}^{n}$. For admissible
configurations $\mathbf{q}\in V$ of the linkage, it is $f_{l}\left( \mathbf{q%
}\right) =\mathbf{I}$. The truncated series expansion (\ref{fTaylor}) thus
defines a $\nu $th-order approximation of the c-space at $\mathbf{q}\in V$ as%
\begin{eqnarray}
V_{\mathbf{q}}^{\nu }{:=} &\{\mathbf{x}\in &{\mathbb{R}}^{n}|\mathrm{d}f_{l,%
\mathbf{q}}\left( \mathbf{x}\right) +\frac{1}{2}\mathrm{d}^{2}f_{l,\mathbf{q}%
}\left( \mathbf{x}\right) +\ldots +\frac{1}{\nu !}\mathrm{d}^{\nu }f_{l,%
\mathbf{q}}\left( \mathbf{x}\right) =\mathbf{0},  \notag \\
&&\ \ \ \ \ \ l=1,\ldots ,\gamma \}.  \label{Vk}
\end{eqnarray}%
The c-space $V$ as a geometric object is independent from the (constraint)
equations used to define it. The $\nu $th-order approximations $V_{\mathbf{q}%
}^{\nu }$ depend on the particular constraints, however. Nevertheless, there
is an order $\kappa $ so that $V_{\mathbf{q}}^{\kappa }$ is a faithful local
approximation of $V$, i.e. there exist a neighborhood $U\left( \mathbf{q}%
\right) $ and an order $\kappa $ such that $V_{\mathbf{q}}^{\kappa }\cap
U\left( \mathbf{q}\right) =V\cap U\left( \mathbf{q}\right) $. 
\vspace{-1ex}%

\section{%
\label{secLocLk}%
Local Approximation of the Set of Constraint Singularities of certain Rank}

\begin{figure*}[t]
\setcounter{equation}{25} 
\begin{equation}
\mathrm{d}^{\nu }m_{%
\bm{\alpha}%
\bm{\beta}%
}=\sum_{\left\vert \mathbf{a}_{k}\right\vert =\nu }\left\vert \mathrm{d}%
^{a_{1}}\mathbf{S}_{%
\bm{\alpha}%
\beta _{1},\mathbf{q}}\left( \mathbf{x}\right) \ \mathrm{d}^{a_{2}}\mathbf{S}%
_{%
\bm{\alpha}%
\beta _{2},\mathbf{q}}\left( \mathbf{x}\right) \ \cdots \ \mathrm{d}^{a_{k}}%
\mathbf{S}_{%
\bm{\alpha}%
\beta _{k},\mathbf{q}}\left( \mathbf{x}\right) \right\vert \frac{\nu !}{%
\mathbf{a}_{k}!n_{1}!\cdots n_{\nu }!}  \label{dm}
\end{equation}%
\setcounter{equation}{15} 
\end{figure*}

The kinematic tangent cone $C_{\mathbf{q}}^{\text{K}}L_{k}$ reveals the
geometry of smooth motions where $\mathbf{J}$ has rank less than $k$. This
fails to provide the c-space geometry at non-tangential intersections of
submanifolds, e.g. cusps. Instead of restricting to motions through a point,
the set $L_{k}$ is now locally approximated.

For $\mathbf{q}\in L_{k}$, all $k$-minors of the constraint Jacobian $%
\mathbf{J}$ in (\ref{J}) vanish, and a $\nu $th-order approximation of $%
L_{k} $ at $\mathbf{q}\in V$ is given by%
\begin{align}
L_{k,\mathbf{q}}^{\nu }=\{\mathbf{x}\in {\mathbb{R}}^{n}|& \mathrm{d}f_{l,%
\mathbf{q}}\left( \mathbf{x}\right) +\frac{1}{2}\mathrm{d}^{2}f_{l,\mathbf{q}%
}\left( \mathbf{x}\right) +\ldots +\frac{1}{\nu !}\mathrm{d}^{\nu }f_{l,%
\mathbf{q}}\left( \mathbf{x}\right) =\mathbf{0},  \notag \\
& \mathrm{d}m_{%
\bm{\alpha}%
\bm{\beta}%
,\mathbf{q}}%
\hspace{-0.5ex}%
\left( \mathbf{x}\right) 
\hspace{-0.5ex}%
+%
\hspace{-0.5ex}%
\frac{1}{2}\mathrm{d}^{2}m_{%
\bm{\alpha}%
\bm{\beta}%
,\mathbf{q}}%
\hspace{-0.5ex}%
\left( \mathbf{x}\right) 
\hspace{-0.5ex}%
+%
\hspace{-0.5ex}%
\ldots 
\hspace{-0.5ex}%
+%
\hspace{-0.5ex}%
\frac{1}{\nu !}\mathrm{d}^{\nu }m_{%
\bm{\alpha}%
\bm{\beta}%
,\mathbf{q}}%
\hspace{-0.5ex}%
\left( \mathbf{x}\right) =\mathbf{0},  \notag \\
& |%
\bm{\alpha}%
|=|%
\bm{\beta}%
|=k\}  \label{Lknu}
\end{align}%
with the differentials of $f_{l}$ and $m_{%
\bm{\alpha}%
\bm{\beta}%
}$, respectively. There is a neighborhood $U\left( \mathbf{q}\right) $ and
order $\kappa $ so that $L_{k,\mathbf{q}}^{\kappa }\cap U\left( \mathbf{q}%
\right) =L_{k}\cap U\left( \mathbf{q}\right) $. All

\section{%
\label{secDiff}%
Differentials and Series Expansions}

In this section, recursive expressions for the differentials appearing in
the definitions of sections \ref{secLocCSpace} and \ref{secLocLk} are given.
The detailed derivations will be presented in a forthcoming publication.
They only require simple vector arithmetic that can be implemented without
the need for symbolic computer algebra systems.%
\vspace{-3ex}%

\subsection{Constraint Mapping}

The above approximation of the c-space $V$ requires the series expansion of
the mapping (\ref{fl}). W.l.o.g. in the following the constraint mapping (%
\ref{f}) of a single loop linkage is considered, for simplicity. The results
are directly applicable to the constraint mapping according to a FC $l$ of a
multi-loop linkage using the predecessor relation $\leq _{l}$ within this FC 
\cite{Robotica2017}. The differentials can be expressed as \cite{Vetter1973}%
\begin{eqnarray}
\mathrm{d}^{k}f_{\mathbf{q}}\left( \mathbf{x}\right) &=&\sum_{\left\vert 
\mathbf{a}\right\vert =k}\frac{k!}{a_{1}!a_{2}!\cdots a_{n}!}\frac{\partial
^{k}f}{\partial q_{1}^{a_{1}}\partial q_{2}^{a_{2}}\ldots \partial
q_{n}^{a_{n}}}x_{1}^{a_{1}}x_{2}^{a_{2}}\cdots x_{n}^{a_{n}}  \notag \\
&=&\sum_{\left\vert \mathbf{a}\right\vert =k}\frac{k!}{\mathbf{a}!}\mathbf{x}%
^{\mathbf{a}}\partial ^{\mathbf{a}}f  \label{dkf1} \\
&=&\sum_{\alpha _{1},\alpha _{2},\ldots ,\alpha _{k}}x_{\alpha
_{1}}x_{\alpha _{2}}\cdots x_{\alpha _{n}}\frac{\partial ^{k}f}{\partial
q_{\alpha _{1}}\partial q_{\alpha _{2}}\ldots \partial q_{\alpha _{k}}}.
\label{dkf2}
\end{eqnarray}%
In (\ref{dkf1}), the multi-index $\mathbf{a}=\left( a_{1},a_{2},\ldots
,a_{n}\right) \in \mathbb{N}^{n}$ is used, where $\mathbf{x}^{\mathbf{a}%
}=x^{a_{1}}x^{a_{2}}\cdots x^{a_{n}}$ is a monomial of degree $k=\left\vert 
\mathbf{a}\right\vert :=a_{1}+a_{2}+\ldots +a_{n}$, and the multiple partial
derivative operator is defined as $\partial ^{\mathbf{a}}=\left( \frac{%
\partial }{\partial q_{1}}\right) ^{a_{1}}\ldots \left( \frac{\partial }{%
\partial q_{n}}\right) ^{a_{n}}$. The multi-factorial is defined as $\mathbf{%
a}!=a_{1}!a_{2}!\cdots a_{n}!$. In the second form (\ref{dkf2}), $\alpha
_{1},\alpha _{2},\ldots ,\alpha _{k}$ are joint indexes running over the
index range $\{1,\ldots ,n\}$.

The differentials can be determined explicitly and algebraically since the
partial derivatives of (\ref{f}) are available in terms of Lie brackets \cite%
{MMTDerivatives}. The non-zero derivatives for $\alpha _{1},\ldots ,a_{\nu
}<i$ are%
\begin{eqnarray}
\frac{\partial ^{\nu }\mathbf{S}_{i}}{\partial q_{\alpha _{1}}\partial
q_{\alpha _{2}}\cdots \partial q_{\alpha _{\nu }}} &=&\left[ \mathbf{S}%
_{\beta _{\nu }},\left[ \mathbf{S}_{\beta _{\nu -1}},\left[ \mathbf{S}%
_{\beta _{\nu -2}},\ldots \left[ \mathbf{S}_{\beta _{1}},\mathbf{S}_{i}%
\right] \ldots \right] \right] \right]  \notag \\
&=&\mathbf{ad}_{\mathbf{S}_{\beta _{\nu }}}\mathbf{ad}_{\mathbf{S}_{\beta
_{\nu -1}}}\mathbf{ad}_{\mathbf{S}_{\beta _{\nu -2}}}\cdots \mathbf{ad}_{%
\mathbf{S}_{\beta _{1}}}\mathbf{S}_{i},  \label{dnS}
\end{eqnarray}%
where $\beta _{\nu }\leq \beta _{\nu -1}\leq \ldots \leq \beta _{2}\leq
\beta _{1}<i$ is the ordered set of indexes $\{\alpha _{1},\ldots ,a_{\nu
}\} $, and $\left[ \mathbf{S}_{j},\mathbf{S}_{k}\right] =\mathbf{ad}_{%
\mathbf{S}_{j}}\mathbf{S}_{k}$ \cite{Selig}. This relation can be
represented using multi-index notation as%
\begin{eqnarray}
\partial ^{\mathbf{a}}\mathbf{S}_{i} &=&\prod\limits_{j=1}^{i-1}\mathbf{ad}_{%
\mathbf{S}_{j}}^{a_{j}}\mathbf{S}_{i},\ \text{if }a_{j}=0,j=i,\ldots ,n
\label{dnS2} \\
&=&\partial ^{\mathbf{a}_{i-1}}\mathbf{S}_{i}
\end{eqnarray}%
and zero otherwise if $a_{j}\neq 0$ for some $j\geq i$.

Implementation of (\ref{dnS}) is straightforward and only involves vector
operations, but the number of operations for its evaluation is rather high.
This can be avoided by means of a recursive formulation.

A straightforward derivation yields the following recursive relations of
differentials of higher degree at $\mathbf{q}\in V$ 
\begin{equation}
\mathrm{d}^{k}f_{\mathbf{q}}\left( \mathbf{x}\right) =h_{\mathbf{q}}^{\left(
k\right) }%
\hspace{-0.5ex}%
\left( \mathbf{x}\right) -\sum_{i=1}^{k-1}\binom{k-1}{i-1}\mathrm{d}^{i}f_{%
\mathbf{q}}\left( \mathbf{x}\right) \mathrm{d}^{k-i}f_{\mathbf{q}}^{-1}%
\hspace{-0.5ex}%
\left( \mathbf{x}\right) ,\text{ for }k>1  \label{dkf}
\end{equation}%
and $\mathrm{d}f_{\mathbf{q}}\left( \mathbf{x}\right) =h_{\mathbf{q}%
}^{\left( 1\right) }%
\hspace{-0.5ex}%
\left( \mathbf{x}\right) $ for $k=1$, where 
\begin{equation}
h_{\mathbf{q}}^{\left( k\right) }%
\hspace{-0.5ex}%
\left( \mathbf{x}\right) :=\sum_{i\leq n}x_{i}\mathrm{d}^{k-1}\widehat{%
\mathbf{S}}_{i,\mathbf{q}}\left( \mathbf{x}\right) .  \label{hk}
\end{equation}%
With slight abuse of notation, $\mathrm{d}f_{\mathbf{q}}^{-1}%
\hspace{-0.5ex}%
\left( \mathbf{x}\right) $ denotes the differential of $f\left( \mathbf{q}%
\right) ^{-1}\in SE\left( 3\right) $. The latter admits the recursive form%
\begin{equation}
\mathrm{d}^{k}f_{\mathbf{q}}^{-1}%
\hspace{-0.5ex}%
\left( \mathbf{x}\right) =-\sum_{i=1}^{k}\binom{k}{i}\mathrm{d}^{i}f_{%
\mathbf{q}}\left( \mathbf{x}\right) \mathrm{d}^{k-i}f_{\mathbf{q}}^{-1}%
\hspace{-0.5ex}%
\left( \mathbf{x}\right) .  \label{dkfinv}
\end{equation}%
Finally, the differentials of the joint screws $\mathbf{S}_{i}$ are given
recursively by%
\begin{equation}
\mathrm{d}^{k}\mathbf{S}_{i,\mathbf{q}}\left( \mathbf{x}\right) =\sum_{j\leq
i}\sum_{r=0}^{k-1}\binom{k-1}{r}[\mathrm{d}^{r}\mathbf{S}_{j,\mathbf{q}%
}\left( \mathbf{x}\right) ,\mathrm{d}^{k-r-1}\mathbf{S}_{i,\mathbf{q}}\left( 
\mathbf{x}\right) ]x_{j},k\geq 1  \label{dSk}
\end{equation}%
with $f^{-1}\left( \mathbf{q}\right) =\mathbf{I}$. In the above expressions,
for $k=0$, it is $\mathrm{d}^{k}\mathbf{S}_{j,\mathbf{q}}=\mathbf{S}%
_{j}\left( \mathbf{q}\right) $ and $\mathrm{d}^{k}f_{\mathbf{q}%
}^{-1}=f^{-1}\left( \mathbf{q}\right) $. It is clear from their definitions
that $\mathrm{d}^{k}f_{\mathbf{q}}\left( \mathbf{x}\right) ,\mathrm{d}^{k}f_{%
\mathbf{q}}^{-1}%
\hspace{-0.5ex}%
\left( \mathbf{x}\right) ,h_{\mathbf{q}}^{\left( k\right) }%
\hspace{-0.5ex}%
\left( \mathbf{x}\right) ,\mathrm{d}^{k}\mathbf{S}_{i,\mathbf{q}}\left( 
\mathbf{x}\right) $ are homogenous polynomials of degree $k$.

With the assumption of remark \ref{rem1}, the origin $\mathbf{q}=\mathbf{0}$
can be chosen as the point of interest. Then (\ref{fTaylor}) becomes the
MacLaurin series, and $\mathbf{S}_{j}\left( \mathbf{0}\right) =\mathbf{Y}%
_{j} $ are the screw coordinates in the reference configuration.

\begin{figure*}[tbh]
\setcounter{equation}{27} 
\begin{align}
&\mathrm{d}f_{\mathbf{q}_{0}}\left( \mathbf{x}\right) =\left( 
\begin{array}{cccc}
0 & -x_{1}-x_{2}-x_{3}-x_{4} & 0 & 0 \\ 
x_{1}+x_{2}+x_{3}+x_{4} & 0 & 0 & -L(x_{1}+2x_{2}+x_{3}) \\ 
0 & 0 & 0 & 0 \\ 
0 & 0 & 0 & 0%
\end{array}%
\right)  \label{df4bar} \\
&\mathrm{d}f_{\mathbf{q}_{0}}\left( \mathbf{x}\right) +\frac{1}{2}\mathrm{d}%
^{2}f_{\mathbf{q}_{0}}\left( \mathbf{x}\right) =\left( 
\begin{array}{cccc}
-\frac{1}{2}(\text{$x_{1}$}+\text{$x_{2}$}+\text{$x_{3}$}+\text{$x_{4}$})^{2}
& -\text{$x_{1}$}-\text{$x_{2}$}-\text{$x_{3}$}-\text{$x_{4}$} & 0 & \frac{1%
}{2}L\left( \text{$x_{1}$}^{2}+4\text{$x_{2}x_{1}$}+2\text{$x_{2}$}^{2}+%
\text{$x_{3}$}^{2}+2(\text{$x_{1}$}+\text{$x_{2}$})\text{$x_{3}$}\right) \\ 
\text{$x_{1}$}+\text{$x_{2}$}+\text{$x_{3}$}+\text{$x_{4}$} & -\frac{1}{2}(%
\text{$x_{1}$}+\text{$x_{2}$}+\text{$x_{3}$}+\text{$x_{4}$})^{2} & 0 & -L(%
\text{$x_{1}$}+2\text{$x_{2}$}+\text{$x_{3}$}) \\ 
0 & 0 & 0 & 0 \\ 
0 & 0 & 0 & 0%
\end{array}%
\right)  \notag \\
&\mathrm{d}f_{\mathbf{q}_{0}}\left( \mathbf{x}\right) +\frac{1}{2}\mathrm{d}%
^{2}f_{\mathbf{q}_{0}}\left( \mathbf{x}\right) +\frac{1}{6}\mathrm{d}^{3}f_{%
\mathbf{q}_{0}}\left( \mathbf{x}\right) =  \notag \\
& {\small \left( 
\begin{array}{cccc}
-\frac{1}{2}(x_{1}+x_{2}+x_{3}+x_{4})^{2} & \frac{1}{6}%
(x_{1}+x_{2}+x_{3}+x_{4})^{3}-x_{1}-x_{2}-x_{3}-x_{4} & 0 & \frac{1}{2}%
L\left(
x_{1}^{2}+4x_{2}x_{1}+2x_{2}^{2}+x_{3}^{2}+2(x_{1}+x_{2})x_{3}\right) 
\vspace{1ex} \\ 
-\frac{1}{6}(x_{1}+x_{2}+x_{3}+x_{4})^{3}+x_{1}+x_{2}+x_{3}+x_{4} & -\frac{1%
}{2}(x_{1}+x_{2}+x_{3}+x_{4})^{2} & 0 & 
\begin{array}{l}
\frac{1}{6}%
L(x_{1}^{3}+6x_{2}x_{1}^{2}+6x_{2}^{2}x_{1}+2x_{2}^{3}+x_{3}^{3}+3(x_{1}+x_{2})x_{3}^{2}
\\ 
\,\,\,\,\,\,\,\,\,\,\,\,\,\,+3(x_{1}+x_{2})^{2}x_{3}-6(x_{1}+2x_{2}+x_{3}))%
\end{array}%
\vspace{1ex} \\ 
0 & 0 & 0 & 0 \\ 
0 & 0 & 0 & 0%
\end{array}%
\right) }  \notag
\end{align}
\setcounter{equation}{26} 
\end{figure*}

\subsection{Jacobian Minors}

Application of the properties of the derivatives of a determinant gives rise
to the general expression (\ref{dm}) for the differentials of the minors,
where the truncated multi-index $\mathbf{a}_{k}=\left( a_{1},a_{2},\ldots
,a_{k}\right) \in \mathbb{N}^{k}$ is defined by removing the last $n-k$
indexes from $\mathbf{a}$, and $n_{i}:=|\{a_{j}|a_{j}=i\}|$ is the number of
times a differential of degree $i$ occurs. The differentials of $\mathbf{S}_{%
\bm{\alpha}%
\beta _{i}}$ are determined by (\ref{dSk}).

\section{Examples%
\label{secExamples}%
}

\subsection{Planar 4-Bar Linkage}

\begin{figure}[b]
{{\includegraphics[width=8.5cm]{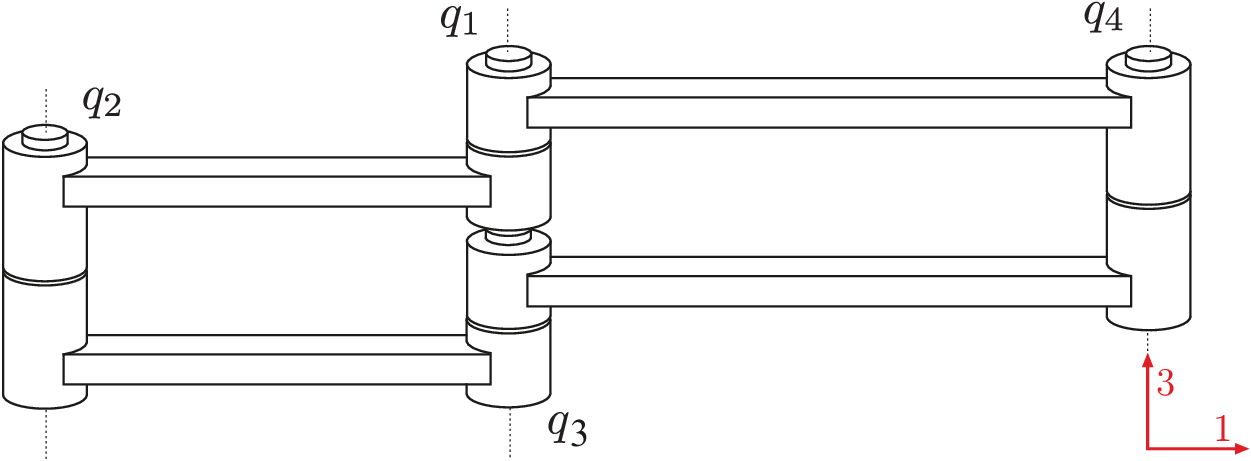}}}
\caption{Reference configuration of the 4-bar linkage (from \protect\cite%
{JMR2016}).}
\label{fig4Bar}
\end{figure}
The first example is the planar 4-bar linkage with four revolute joints,
i.e. $V\subset {\mathbb{T}}^{4}$. In the (singular) reference configuration $%
\mathbf{q}_{0}=\mathbf{0}$, shown in fig. \ref{fig4Bar}, the $n=4$ joint
screw coordinates are%
\begin{eqnarray*}
\mathbf{Y}_{1} &=&\mathbf{Y}_{3}=(0,0,1,0,L,0)^{T},\mathbf{Y}%
_{2}=(0,0,1,0,2L,0)^{T} \\
\mathbf{Y}_{4} &=&(0,0,1,0,0,0)^{T}.
\end{eqnarray*}

\paragraph{Analysis of $V$:}

The kinematic tangent cone to $V$ was shown to be \cite{JMR2016}%
\begin{equation*}
C_{\mathbf{q}_{0}}^{\text{K}}V=\mathbf{V}(x_{1}+x_{3},x_{2},x_{4})\cup 
\mathbf{V}(x_{1}+x_{2},x_{3}+x_{2},x_{4}-x_{2}).
\end{equation*}%
That is, the linkage can perform smooth 1-DOF motions through $\mathbf{q}%
_{0} $, and two motion branches intersect at that point.

It is known that the configuration $\mathbf{q}_{0}$ is a bifurcation of $V$,
and $C_{\mathbf{q}_{0}}^{\text{K}}V$ is in fact the union of the two tangent
spaces to the manifolds intersecting at $\mathbf{q}_{0}$. This cannot be
inferred from the kinematic tangent cone, however, but from a local
approximation of $V$. To this end, the $i$th-order series expansions of the
constraint mapping are determined using (\ref{dkf}). The expansions up to
3rd-order are shown in (\ref{df4bar}). They, and the higher-order
expansions, yield the approximations (\ref{Vk}) as $V_{\mathbf{q}%
_{0}}^{i}=C_{\mathbf{q}_{0}}^{\text{K}}V,i\geq 2$. This confirms that $%
\mathbf{q}_{0}$ is a bifurcation point of $V$ and thus a c-space
singularity. The local DOF is $\delta _{\mathrm{loc}}\left( \mathbf{q}%
_{0}\right) =\dim V_{\mathbf{q}_{0}}^{2}=1$.

\paragraph{Analysis of $L_{i}$:}

It remains to investigate the differential DOF at configurations near $%
\mathbf{q}_{0}$. The rank of the Jacobian $\mathbf{J}\left( \mathbf{q}%
_{0}\right) =\left( \mathbf{Y}_{1}\ \mathbf{Y}_{2}\ \mathbf{Y}_{3}\ \mathbf{Y%
}_{4}\right) $ in the reference configuration is $\mathrm{rank}~\mathbf{J}%
\left( \mathbf{q}_{0}\right) =2$, i.e. the differential DOF is $\delta _{%
\text{diff}}\left( \mathbf{q}_{0}\right) =n-2$, so that $\mathbf{q}_{0}\in
L_{3}$. For this planar linkage the maximal rank of $\mathbf{J}$ is 3 so
that trivially $\mathbf{q}_{0}\in L_{i},i>3$. The kinematic tangent cone to $%
L_{3}$ was reported in \cite{JMR2018} to be $C_{\mathbf{q}_{0}}^{\text{K}%
}L_{3}=\{\mathbf{0}\}$. From this analysis it can be concluded that for any
smooth curve through $\mathbf{q}_{0}$ the differential DOF is locally $%
\delta _{\mathrm{diff}}=1$ but increases to 2 at $\mathbf{q}_{0}$.

The local approximation of $L_{3}$ requires the differentials of the
2-minors of $\mathbf{J}$ at $\mathbf{q}_{0}$. The first- and third-order
differentials are%
\begin{eqnarray*}
\{\mathrm{d}m_{%
\bm{\alpha}%
\bm{\beta}%
,\mathbf{q}_{0}}%
\hspace{-0.5ex}%
\left( \mathbf{x}\right) \} &=&\left\{ -L^{2}\text{$x_{2}$},-L^{2}(2\text{$%
x_{2}$}+\text{$x_{3}$}),-L^{2}\text{$x_{2}$},L^{2}\text{$x_{3}$}\right\} \\
\{\mathrm{d}^{3}m_{%
\bm{\alpha}%
\bm{\beta}%
,\mathbf{q}_{0}}%
\hspace{-0.5ex}%
\left( \mathbf{x}\right) \} &=&\{L^{2}\text{$x_{2}$}^{3},L^{2}\left( 2\text{$%
x_{2}$}^{3}+3\text{$x_{2}$}^{2}\text{$x_{3}$}+3\text{$x_{2}x_{3}$}^{2}+\text{%
$x_{3}$}^{3}\right) , \\
&&\ \ -L^{2}\text{$x_{3}$}^{3},L^{2}\text{$x_{2}$}\left( \text{$x_{2}$}^{2}+3%
\text{$x_{2}x_{3}$}+3\text{$x_{3}$}^{2}\right) \},
\end{eqnarray*}%
whereas $\mathrm{d}^{2}m_{%
\bm{\alpha}%
\bm{\beta}%
,\mathbf{q}_{0}}%
\hspace{-0.5ex}%
\left( \mathbf{x}\right) =0$, for $|%
\bm{\alpha}%
|=|%
\bm{\beta}%
|=2$. Therewith the solution sets (\ref{Lknu}) are $L_{3,\mathbf{q}%
_{0}}^{i}=C_{\mathbf{q}_{0}}^{\text{K}}L_{3}=\{\mathbf{0}\},i>1$.
Consequently, all points in a neighborhood of $\mathbf{q}_{0}=\mathbf{0}$
have rank 3, except the point $\mathbf{q}_{0}$. The set of points of rank
less than 3 (and thus of rank 2) is the zero-dimensional manifold consisting
of the point $\{\mathbf{q}_{0}\}$.%
\vspace{-2ex}%

\begin{figure}[b]
{{a)\centerline{%
\includegraphics[width=8.5cm]{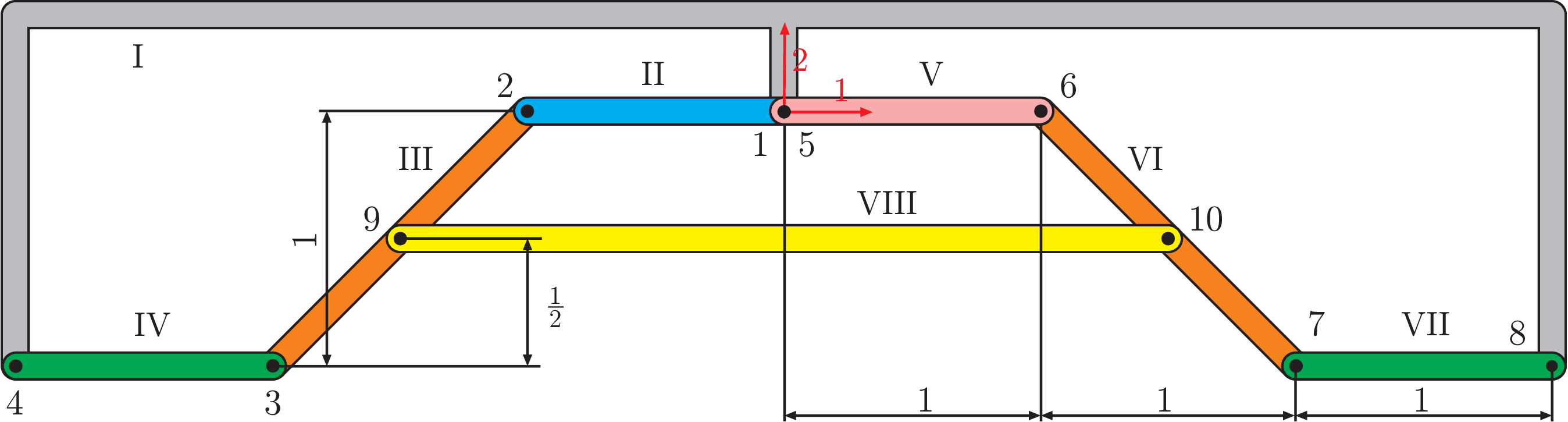}} \vspace{%
1ex} }}
\par
b)\centerline{%
\includegraphics[width=5.5cm]{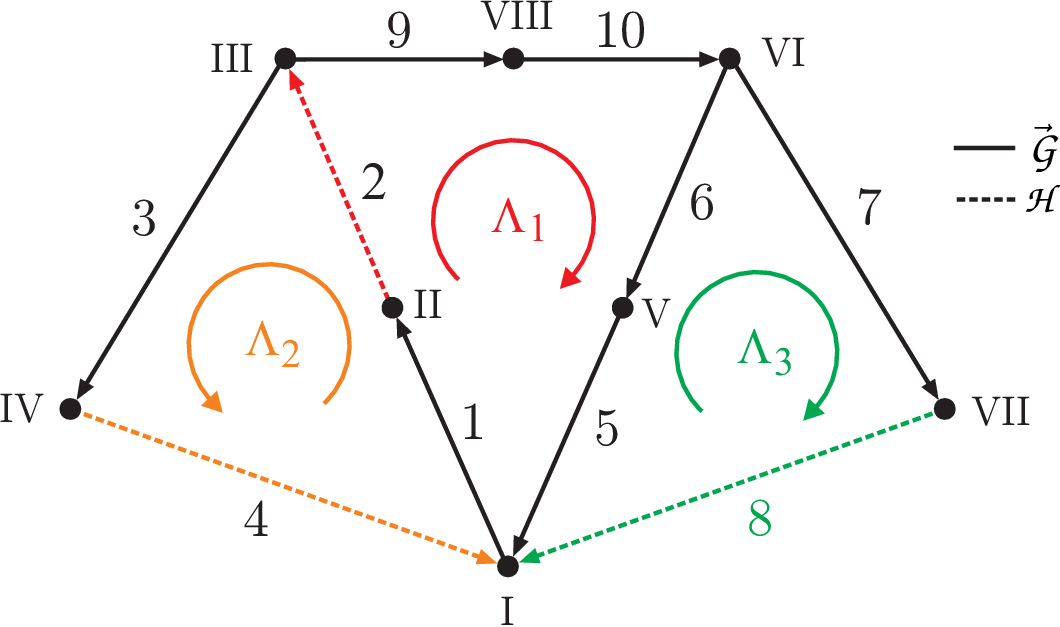}}
\caption{a) Planar linkage, constructed by connecting two Watt linkages, in
its reference configuration $\mathbf{q}_{0}$ (from \protect\cite{JMR2018}).
b) Topological graph and fundamental cycles.}
\label{figDoubleWatt}
\end{figure}

\subsection{Planar 3-Loop Linkage in a Cusp Singularity}

The 1 DOF planar linkage shown in fig. \ref{figDoubleWatt} was discussed in 
\cite{ConnellyServatius1994} as an example where the c-space exhibits a cusp
singularity. It is formed by a mirror-symmetric arrangement of two
straight-line-generating Watt linkages that are connected by link VIII. The
attachment point of joint 9 moves on the coupler point of a Watt linkages
formed by links I-IV, and the attachment point of joint 10 moves on the
coupler curve of a mirrored version of this Watt linkages formed by links I,
V-VII. Therefore, both points can instantaneously only move in vertical
direction. Furthermore, in the shown configuration, both points are at an
inflection point of the respective coupler curve. Since the sitance of these
points is fixed by link VIII, link III or VI can instantaneously not move
upward. There are two possible downward motions, which lead to the cusp of
the c-space. A detailed description can be found in \cite%
{ConnellyServatius1994,PabloMMT2019}.

It comprises $n=10$ revolute joints with parallel axes, i.e. $V\subset {%
\mathbb{T}}^{10}$. In the reference configuration in fig. \ref{figDoubleWatt}
the joint screw coordinates, w.r.t. the shown reference frame, are%
\begin{eqnarray*}
\mathbf{Y}_{1} &=&\mathbf{Y}_{5}=(0,0,1,0,0,0)^{T},\mathbf{Y}_{2}=\left(
0,0,1,0,1,0\right) ^{T} \\
\mathbf{Y}_{3} &=&(0,0,1,-1,2,0)^{T},\mathbf{Y}_{4}=(0,0,1,-1,3,0)^{T} \\
\mathbf{Y}_{6} &=&(0,0,1,0,-1,0)^{T},\mathbf{Y}_{7}=(0,0,1,-1,-2,0)^{T} \\
\mathbf{Y}_{8} &=&(0,0,1,-1,-3,0)^{T},\mathbf{Y}_{9}=(0,0,1,-1/2,3/2,0)^{T}
\\
\mathbf{Y}_{10} &=&(0,0,1,-1/2,-3/2,0)^{T}.%
\vspace{-2ex}%
\end{eqnarray*}%
The linkage possesses $\gamma =3$ FCs numbered with 1,2,3 (instead of using
indexes of the co-tree edge as in \cite{Robotica2017,JMR2018}). The loop
constraint mappings are%
\begin{eqnarray*}
f_{1}\left( \mathbf{q}\right) &=&\exp \left( q_{1}\mathbf{Y}_{1}\right) \exp
\left( q_{2}\mathbf{Y}_{2}\right) \exp \left( q_{9}\mathbf{Y}_{9}\right)
\cdot \\
&&\exp \left( q_{10}\mathbf{Y}_{10}\right) \exp \left( q_{6}\mathbf{Y}%
_{6}\right) \exp \left( q_{5}\mathbf{Y}_{5}\right) \\
f_{2}\left( \mathbf{q}\right) &=&\exp \left( q_{1}\mathbf{Y}_{1}\right) \exp
\left( q_{2}\mathbf{Y}_{2}\right) \exp \left( q_{3}\mathbf{Y}_{3}\right)
\exp \left( q_{4}\mathbf{Y}_{4}\right) \\
f_{3}\left( \mathbf{q}\right) &=&\exp \left( -q_{5}\mathbf{Y}_{5}\right)
\exp \left( -q_{6}\mathbf{Y}_{6}\right) \exp \left( q_{7}\mathbf{Y}%
_{7}\right) \exp \left( q_{8}\mathbf{Y}_{8}\right) .
\end{eqnarray*}%
The corresponding constraint Jacobians are\setcounter{equation}{28} 
\begin{eqnarray}
\mathbf{J}_{1} &=&%
\Big%
(\mathbf{S}_{1},\mathbf{S}_{2},\mathbf{0},\mathbf{0},\mathbf{S}_{5},\mathbf{S%
}_{6},\mathbf{0},\mathbf{0},\mathbf{S}_{9},\mathbf{S}_{10}%
\Big%
)  \notag \\
\mathbf{J}_{2} &=&%
\Big%
(\mathbf{S}_{1},\mathbf{S}_{2},\mathbf{S}_{3},\mathbf{S}_{4},\mathbf{0},%
\mathbf{0},\mathbf{0},\mathbf{0},\mathbf{0},\mathbf{0}%
\Big%
)  \notag \\
\mathbf{J}_{3} &=&%
\Big%
(\mathbf{0},\mathbf{0},\mathbf{0},\mathbf{0},-\mathbf{S}_{5},-\mathbf{S}_{6},%
\mathbf{S}_{7},\mathbf{S}_{8},\mathbf{0},\mathbf{0}%
\Big%
).%
\vspace{-4ex}%
\end{eqnarray}%
They form the overall $18\times 10$ constraint Jacobian (\ref{J}). In the
reference configuration, where $\mathbf{S}_{i}\left( \mathbf{q}_{0}\right) =%
\mathbf{Y}_{i}$, it is $r=\mathrm{rank}~\mathbf{J}\left( \mathbf{q}%
_{0}\right) =8$, and hence $\delta _{\text{diff}}\left( \mathbf{q}%
_{0}\right) =2$.

\paragraph{Analysis of $V$:}

It was shown in \cite{JMR2018} that the kinematic tangent cone to the
c-space is ${C_{\mathbf{q}_{0}}^{\text{K}}V}=\{\mathbf{0}\}$. That is, there
is no smooth curve through $\mathbf{q}_{0}$. Yet the linkage is known to be
finitely mobile with 1 DOF, i.e. the c-space $V$ is a curve in ${\mathbb{T}}%
^{10}$. This can only be revealed by analyzing the geometry of $V$.

\begin{figure}[t!]
\begin{center}
a)%
\includegraphics[width=6cm]{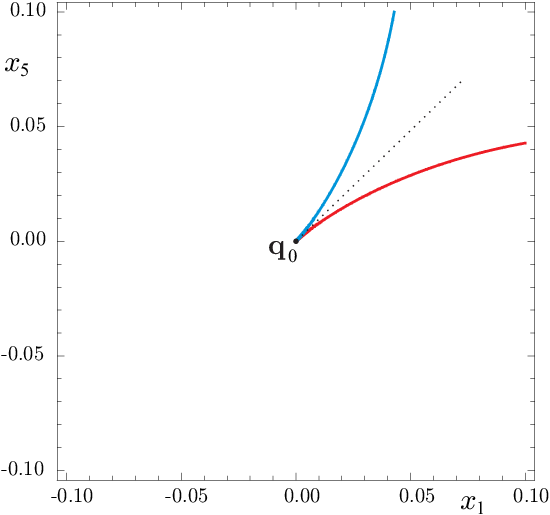}
\vspace{2ex}\\[0pt]
b)%
\includegraphics[width=6cm]{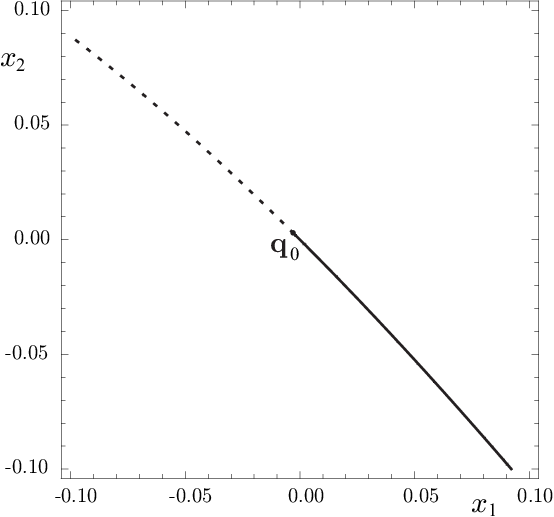}
\end{center}
\caption{2nd-order approximation $V_{\mathbf{q}_{0}}^{2}$ of the c-space of
the linkage in fig. \protect\ref{figDoubleWatt}: a) $x_{1}$-$x_{5}$ section,
b) $x_{1}$-$x_{2}$ section.}
\label{figDoubleWattCSpace1}
\end{figure}

The first-order approximation (\ref{Vk}) is simply the 2-dimensional vector
space 
\begin{eqnarray}
V_{\mathbf{q}_{0}}^{1} &=&\ker \mathbf{J}\left( \mathbf{q}_{0}\right) ={K_{%
\mathbf{q}_{0}}^{1}}  \notag \\
&=&\mathbf{V}%
\big%
(x_{1}+x_{2},x_{1}+x_{3},x_{1}-x_{4},x_{5}+x_{6},x_{5}-x_{7}, \\
&&\ \ \ \ x_{5}+x_{8},x_{1}-x_{5}-3x_{9},x_{1}+3x_{10}-x_{5}%
\big%
).  \notag
\end{eqnarray}%
\begin{figure}[th]
\begin{center}
a)%
\includegraphics[width=8.2cm]{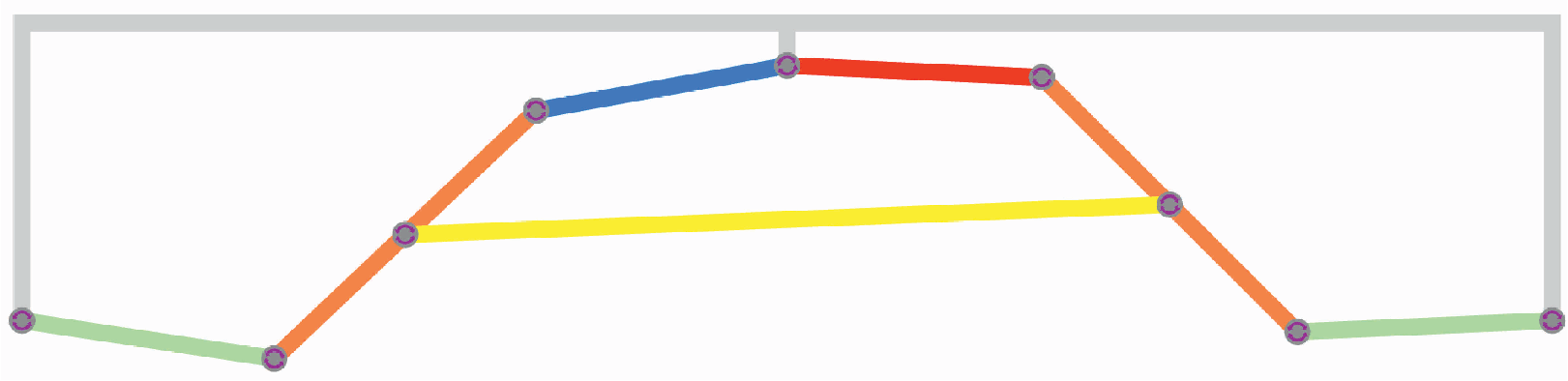}
\vspace{2ex}\\[0pt]
b)%
\includegraphics[width=8.2cm]{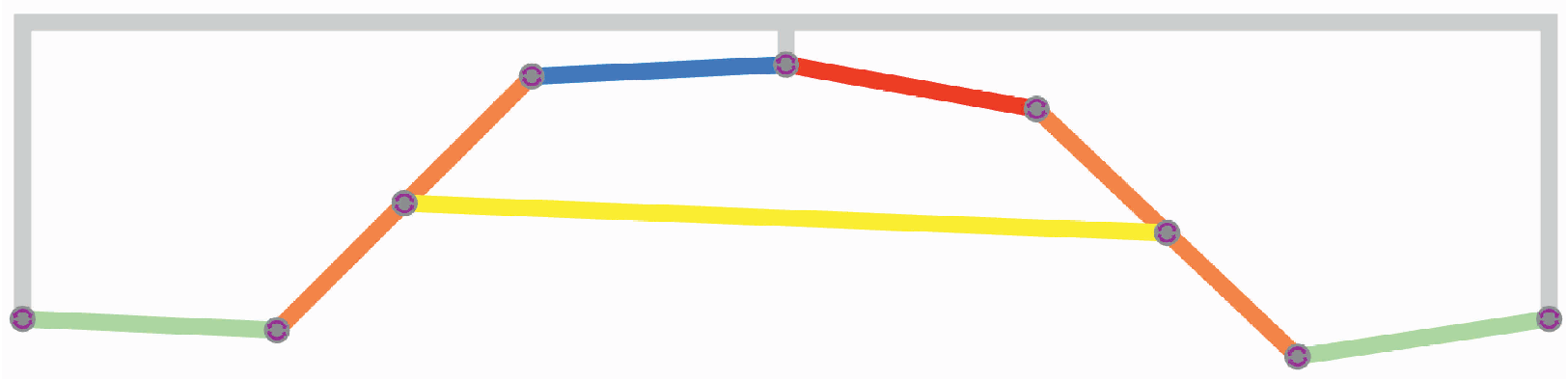}
\end{center}
\caption{Two nearby configurations the linkage can reach from the cusp
singulaity $\mathbf{q}_{0}$, where a) corresponds to the red part, and b) to
the blue part of the curve in Fig. \protect\ref{figDoubleWattCSpace1}a).}
\label{figDoubleWattConfigs}
\end{figure}
\begin{figure*}[tbh]
{\small 
\begin{align}
& \mathrm{d}f_{2,\mathbf{q}_{0}}\left( \mathbf{x}\right) +\frac{1}{2}\mathrm{%
d}^{2}f_{2,\mathbf{q}_{0}}\left( \mathbf{x}\right) =  \notag \\
& \left( 
\begin{array}{cccc}
-\frac{1}{2}({x_{1}}+{x_{10}}+{x_{2}}+{x_{5}}+{x_{6}}+{x_{9}})^{2} & -{x_{1}}%
-{x_{10}}-{x_{2}}-{x_{5}}-{x_{6}}-{x_{9}} & 0 & 
\begin{array}{l}
\frac{1}{4}(3{x_{10}}^{2}+(6{x_{1}}+6{x_{2}}+4{x_{6}}+6{x_{9}}-2){x_{10}}-2{%
x_{2}}^{2}+2{x_{6}}^{2} \\ 
-3{x_{9}}^{2}+4{x_{2}}{x_{6}}-6{x_{2}}{x_{9}}+4{x_{6}}{x_{9}}-2{x_{9}}+{x_{1}%
}(-4{x_{2}}+4{x_{6}}-6{x_{9}}))%
\end{array}
\\ 
{x_{1}}+{x_{10}}+{x_{2}}+{x_{5}}+{x_{6}}+{x_{9}} & -\frac{1}{2}({x_{1}}+{%
x_{10}}+{x_{2}}+{x_{5}}+{x_{6}}+{x_{9}})^{2} & 0 & -\frac{3{x_{10}}}{2}+{%
x_{2}}-{x_{6}}+\frac{3{x_{9}}}{2}-\frac{1}{4}({x_{10}}+{x_{9}})(2{x_{1}}+{%
x_{10}}+2{x_{2}}+{x_{9}}) \\ 
0 & 0 & 0 & 0 \\ 
0 & 0 & 0 & 0%
\end{array}%
\right)  \notag \\
& \mathrm{d}f_{4,\mathbf{q}_{0}}\left( \mathbf{x}\right) +\frac{1}{2}\mathrm{%
d}^{2}f_{4,\mathbf{q}_{0}}\left( \mathbf{x}\right) =  \notag \\
& \left( 
\begin{array}{cccc}
-\frac{1}{2}({x_{1}}+{x_{2}}+{x_{3}}+{x_{4}})^{2} & -{x_{1}}-{x_{2}}-{x_{3}}-%
{x_{4}} & 0 & \frac{1}{2}\left( -{x_{2}}^{2}-2{x_{1}}{x_{2}}-4{x_{3}}{x_{2}}%
-2{x_{3}}^{2}-3{x_{4}}^{2}-4{x_{1}}{x_{3}}-2{x_{3}}-6({x_{1}}+{x_{2}}+{x_{3}}%
){x_{4}}-2{x_{4}}\right) \\ 
{x_{1}}+{x_{2}}+{x_{3}}+{x_{4}} & -\frac{1}{2}({x_{1}}+{x_{2}}+{x_{3}}+{x_{4}%
})^{2} & 0 & {x_{2}}+2{x_{3}}+3{x_{4}}-\frac{1}{2}({x_{3}}+{x_{4}})(2{x_{1}}%
+2{x_{2}}+{x_{3}}+{x_{4}}) \\ 
0 & 0 & 0 & 0 \\ 
0 & 0 & 0 & 0%
\end{array}%
\right)  \notag \\
& \mathrm{d}f_{8,\mathbf{q}_{0}}\left( \mathbf{x}\right) +\frac{1}{2}\mathrm{%
d}^{2}f_{8,\mathbf{q}_{0}}\left( \mathbf{x}\right) =  \label{df2Watt} \\
& \left( 
\begin{array}{cccc}
-\frac{1}{2}({x_{5}}+{x_{6}}-{x_{7}}-{x_{8}})^{2} & {x_{5}}+{x_{6}}-{x_{7}}-{%
x_{8}} & 0 & {x_{5}}({x_{6}}-2{x_{7}}-3{x_{8}})+\frac{1}{2}\left( {x_{6}}%
^{2}-4{x_{7}}{x_{6}}-6{x_{8}}{x_{6}}+2{x_{7}}^{2}+3{x_{8}}^{2}-2{x_{7}}+6{%
x_{7}}{x_{8}}-2{x_{8}}\right) \\ 
-{x_{5}}-{x_{6}}+{x_{7}}+{x_{8}} & -\frac{1}{2}({x_{5}}+{x_{6}}-{x_{7}}-{%
x_{8}})^{2} & 0 & {x_{6}}-2{x_{7}}-3{x_{8}}+\frac{1}{2}(2{x_{5}}+2{x_{6}}-{%
x_{7}}-{x_{8}})({x_{7}}+{x_{8}}) \\ 
0 & 0 & 0 & 0 \\ 
0 & 0 & 0 & 0%
\end{array}%
\right)  \notag
\end{align}%
} 
\end{figure*}
According to (\ref{Vk}), the second-order series expansion of the constraint
mappings, shown in (\ref{df2Watt}), yields the second-order approximation (%
\ref{V2}). This is a 1-dimensional real algebraic variety. It turns out that
the third- and higher-order approximations are all 1-dimensional, and the
local DOF is thus $\delta _{\mathrm{loc}}\left( \mathbf{q}_{0}\right) =\dim
V_{\mathbf{q}_{0}}^{2}=1$. In order to visualize this 1-dimensional variety,
a bivariate polynomial is constructed from (\ref{df2Watt}) in terms of $%
x_{1} $ and $x_{5}$ by eliminating $x_{2},\ldots ,x_{4},x_{6},\ldots ,x_{10}$%
. This is a polynomial of degree 32, which is too long to be shown here. It
defines the $x_{1}$-$x_{5}$ section shown in fig. \ref{figDoubleWattCSpace1}%
a), which reveals that $\mathbf{q}_{0}$ is a cusp singularity. Fig. \ref%
{figDoubleWattConfigs} shows two nearby configurations the linkage can reach
from the cusp $\mathbf{q}_{0}$.

The bivariate polynomials in terms of $x_{1}$ and the other variables are of
degree 4, respectively. For instance,%
\begin{eqnarray*}
x_{1}^{4}+4x_{1}^{3}x_{2}+6x_{1}^{2}x_{2}^{2}+4x_{1}x_{2}^{3}+x_{2}^{4}-8x_{1}^{3}-20x_{1}^{2}x_{2} &&
\\
-16x_{1}x_{2}^{2}-4x_{2}^{3}+24x_{1}^{2}+24x_{1}x_{2}+8x_{2}^{2}+8x_{1}+8x_{2} &=&0
\end{eqnarray*}%
defines the projection of $V_{\mathbf{q}_{0}}^{2}$ to the $x_{1}$-$x_{2}$
coordinate plane shown in fig. \ref{figDoubleWattCSpace1}b). It is important
to notice that the dashed part of the curve is not connected to the cusp,
but rather appears due to the projection, i.e. due to the elimination of the
remaining variables. 
\begin{figure*}[th]
\begin{equation}
\begin{array}{ll}
V_{\mathbf{q}_{0}}^{2}\hspace{-0.6ex}=\hspace{-0.6ex}\mathbf{V}\big(\hspace{%
-0.9ex} & {x_{5}}+{x_{6}}-{x_{7}}-{x_{8}},-\frac{1}{2}({x_{5}}+{x_{6}}-{x_{7}%
}-{x_{8}})^{2},-{x_{5}}-{x_{6}}+{x_{7}}+{x_{8}},\frac{1}{2}({x_{7}}+{x_{8}}%
)(2{x_{5}}+2{x_{6}}-{x_{7}}-{x_{8}})+{x_{6}}-2{x_{7}}-3{x_{8}}, \\ 
& {x_{5}}({x_{6}}-2{x_{7}}-3{x_{8}})+\frac{1}{2}\left( {x_{6}}^{2}-4{x_{6}}{%
x_{7}}-6{x_{6}}{x_{8}}+2{x_{7}}^{2}+6{x_{7}}{x_{8}}-2{x_{7}}+3{x_{8}}^{2}-2{%
x_{8}}\right) ,\frac{1}{4}(-{x_{10}}-4{x_{3}}-4{x_{4}}-{x_{9}}), \\ 
& \frac{1}{2}(-3{x_{1}}-{x_{10}}-3{x_{2}}-2{x_{3}}-2{x_{4}}-{x_{5}}-{x_{6}}-{%
x_{9}}),{x_{1}}+{x_{10}}+{x_{2}}+{x_{5}}+{x_{6}}+{x_{9}},-\frac{1}{2}({x_{1}}%
+{x_{10}}+{x_{2}}+{x_{5}}+{x_{6}}+{x_{9}})^{2}, \\ 
& \frac{1}{4}(-3{x_{10}}+6{x_{2}}+8{x_{3}}+12{x_{4}}-2{x_{6}}+3{x_{9}}),{%
x_{2}}+\frac{3{x_{9}}}{2}-\frac{1}{4}({x_{10}}+{x_{9}})(2{x_{1}}+{x_{10}}+2{%
x_{2}}+{x_{9}})-\frac{3{x_{10}}}{2}-{x_{6}}, \\ 
& \frac{1}{4}\left( {x_{10}}(6{x_{1}}+6{x_{2}}+4{x_{6}}+6{x_{9}}-2)+{x_{1}}%
(-4{x_{2}}+4{x_{6}}-6{x_{9}})+3{x_{10}}^{2}-2{x_{2}}^{2}+4{x_{2}}{x_{6}}-6{%
x_{2}}{x_{9}}+2{x_{6}}^{2}+4{x_{6}}{x_{9}}-3{x_{9}}^{2}-2{x_{9}}\right) , \\ 
& \frac{1}{2}({x_{1}}+{x_{10}}+{x_{2}}+{x_{5}}+{x_{6}}+{x_{9}})+{x_{1}}+{%
x_{2}}+{x_{3}}+{x_{4}},-{x_{1}}-{x_{10}}-{x_{2}}-{x_{5}}-{x_{6}}-{x_{9}}\big)%
\end{array}
\label{V2}
\end{equation}%
\end{figure*}

\begin{remark}
Clearly, the notion of a tangent space or tangent cone to $V$ is
questionable at this point. One relevant definition, for instance, is that
of the geometric tangent cone defined as the limits of all secant lines
originating from $\mathbf{q}$ that pass through a sequence of points that
converges to $\mathbf{q}$ \cite{OSheaWilson,Whitney1965}. This would be the
doted line shown in fig. \ref{figDoubleWattCSpace1}a), which shows that the
linkage is at a dead point. It does, however, not indicate that there are
two branches. Moreover, any definition of tangent cone to the c-space (as
real set) may have a different dimension than $V$ itself, and hence not
reveal the correct DOF.
\end{remark}

\begin{remark}
While the determination of the polynomials defining the second-order
approximation is computationally simple with the recursive relations
presented in section \ref{secDiff}, their interpretation and the actual
analysis of the second-order algebraic variety $V_{\mathbf{q}_{0}}^{2}$ in (%
\ref{V2}) is non-trivial. The polynomial was determined with the software 
\texttt{Maple} on a high-performance computing facility.
\end{remark}

\paragraph{Analysis of $L_{i}$:}

In the reference configuration it is $\mathrm{rank}~\mathbf{J}\left( \mathbf{%
q}_{0}\right) =8$, and since the maximal rank of $\mathbf{J}$ is 9
(recalling that the linkage consists of three planar loops), it is $\mathbf{q%
}_{0}\in L_{9}$. Because any $C_{\mathbf{q}_{0}}^{\text{K}}L_{i}$ is a
restriction of ${C_{\mathbf{q}_{0}}^{\text{K}}V}=\{\mathbf{0}\}$ it is $C_{%
\mathbf{q}_{0}}^{\text{K}}L_{9}=\{\mathbf{0}\}$ (there is no smooth curve
through $\mathbf{q}_{0}$).

The second-order local approximation of the variety $L_{9}$ turns out to be $%
L_{9,\mathbf{q}_{0}}^{2}=\{\mathbf{0}\}$. Therefore, $\mathrm{rank}~\mathbf{J%
}=9$ in any sufficiently small neighborhood $U\left( \mathbf{q}_{0}\right)
\cap V$ except at $\mathbf{q}_{0}$, where the rank drops to 8. (Notice that $%
L_{8}\cap U\left( \mathbf{q}_{0}\right) =\emptyset $ since $\mathrm{rank}~%
\mathbf{J}$ can only increase in a sufficiently small neighborhood). Thus, $%
\mathbf{q}_{0}$ is a c-space singularity, and hence a kinematic singularity.

This example shows that the mobility/singularity analysis in general
requires resorting to the most obvious tool, namely the approximation of the
solution set $V$. However, it also shows that the final analysis amounts to
investigating an algebraic variety $V_{\mathbf{q}_{0}}^{\kappa }$, which
gets complicated very easily.%
\vspace{-3ex}%

\section{Conclusions and Open Issues}

Investigating the finite motion of a linkage strictly requires analyzing the
c-space geometry. This was addressed in this paper by means of a local
approximation. To this end, a recursive algebraic formulation of the
higher-order series expansion of the loop constraints is provided, and used
to define higher-order local approximations of the c-space and thus the
mobility. Additionally, a computationally efficient formulation of the
series expansion of the minors of the constraint Jacobian is presented and
thereupon a local approximation of the set of points of certain rank is
defined. This allows for analyzing configurations where a linkage does not
possess smooth motions through that configuration. It should be noticed that
such situations were almost exclusively ignored in the literature since
singularities of most linkages are bifurcations were smooth manifolds
intersect (even for all reported kinematotropic linkages --a kinematotropic
linkage that exhibits non-smooth transitions was reported in \cite%
{PabloMMT2019}). The presented approach complements the previously presented
formulation for approximating finite motions through a given point.

In summary all defining polynomial systems implicitly defining a local
description of 1) the set of smooth motions, 2) the set of smooth motions
with certain rank (including singularities), 3) the c-space, and 4) the set
of points of certain rank are given explicitly in terms of simple algebraic
operations of joint screw coordinates. Thus, a conceptual and algorithmic
framework providing all defining equations for an exhaustive local mobility
analysis of linkages is now available.

Still the critical point of any local approximation is that it leads to an
approximating algebraic variety defined implicitly by a polynomial system.
Investigating the latter can become difficult, and may necessitate
application of methods from computational algebraic geometry. This may seem
to put the usefulness of the local approximation as whole into perspective
since the c-space itself could be expressed as an algebraic variety in the
first place. However, the order of the approximating algebraic variety is
low in contrast to the algebraic c-space variety. With the formulation
presented in this paper along with \cite{JMR2016,JMR2018}, the defining
equations for a local approximation can always be determined. The best means
to treat these equations remains to be investigated.

The higher-order series expansion can also be used for computing approximate
solutions of loop closure constraints, e.g. for a specified motion. Initial
results are reported in \cite{deJong2018}. This topic has also been
addressed in \cite{Milenkovic2012}.

\section*{Acknowledgement}

This work has been supported by the LCM-K2 Center within the framework of
the Austrian COMET-K2 program.

\end{document}